\documentclass[12pt]{article}
\usepackage{amssymb,latexsym,amsmath,euscript}\usepackage{amsfonts}
\usepackage[cp1251]{inputenc}\usepackage[epsf]{}\usepackage[english,russian]{babel}\usepackage{graphicx}
\newtheorem{theorem}{Theorem}
\newtheorem{lemma}{Lemma}[section]
\newtheorem{proposition}{Proposition}
\newtheorem{corollary}{Corollary}[section]
\begin{document}

\title{On existence of expanding attractors with different dimensions
\footnote{2000{\it Mathematics Subject Classification}. Primary 37D15; Secondary 58C30}
}

\author{V.~Medvedev$^{1}$\and E.~Zhuzhoma$^{1}$}
\date{}
\maketitle

{\small $^{1}$ National Research University Higher School of Economics, 25/12 Bolshaya Pecherskaya, 603005, \\ Nizhny Novgorod, Russia}

\renewcommand{\abstractname}{Absrtact}
\renewcommand{\refname}{Bibliography}
\renewcommand{\figurename}{Fig.}
%%%%%%%%%%%%%%%%%%%%%%%%%%%%%%%%%%%%

\begin{abstract}
We prove that $n$-sphere $\mathbb{S}^n$, $n\geq 2$, admits structurally stable diffeomorphisms $\mathbb{S}^n\to\mathbb{S}^n$ with non-orientable expanding attractors of any topological dimension $d\in\{1,\ldots,[\frac{n}{2}]\}$ where $[x]$ is an integer part of $x$. One proves that $n$-torus $\mathbb{T}^n$, $n\geq 2$, admits structurally stable diffeomorphisms $\mathbb{T}^n\to\mathbb{T}^n$ with orientable expanding attractors of any topological dimension $1\leq q\leq n-1$. We also prove that given any closed $n$-manifold $M^n$, $n\geq 2$, and any $d\in\{1,\ldots,[\frac{n}{2}]\}$, there is an axiom A diffeomorphism $f: M^n\to M^n$ with a $d$-dimensional non-orientable expanding attractor. Similar statements hold for axiom A flows.
\end{abstract}

\section*{Introduction}

Below, $M^n$ is a closed smooth connected $n$-manifold, $n\geq 2$. Let $f: M^n\to M^n$ be a diffeomorphism and $\Lambda$ an invariant set of $f$, i.e. $f(\Lambda)=\Lambda$. This set is called an \textit{attractor} provided there is a neighborhood $U(\Lambda)=U$ of $\Lambda$ such that $\overline{f(U)}\subset U$ and $\cap_{i\geq 0}f^i(U)=\Lambda$ where $\overline{N}$ means the topological closure of $N$.
For applications, the most interesting invariant sets of a dynamical system are attractors and repellers. The famous attractors in hyperbolic dynamical systems are Smale's solenoid, DA-attractor, and Plykin's attractor (see the books \cite{ABZ,GrinesZh-book-2021,Robinson-book-99}). All this attractors are expanding ones and basic sets of axiom A diffeomorphisms. Recall that a diffeomorp\-hism $f: M^n\to M^n$ satisfies an axiom A (in short, $f$ is an A-diffeomorphism) provided a non-wandering set $NW(f)$ of $f$ is hyperbolic, and $NW(f)$ is the topological closure of periodic orbits of $f$ \cite{Smale67}. According Smale's Spectral Decomposition Theorem, the non-wandering set $NW(f)$ splits into pairwise disjoint invariant closed and transitive sets called \textit{basic sets}.

The hyperbolic structure implies the existence of stable $W^s(x)$ and unstable $W^u(x)$ manifolds at any points $x\in NW(f)$. Following \cite{Williams1974} we call a basic set $\Lambda$ an \textit{expanding attractor} provided $\Lambda$ is an attractor and $\dim\Lambda=\dim W^u(x)$ for every point $x\in\Lambda$. Williams \cite{Williams1974} proved that an expanding attractor is locally homeomorphic to the product of Cantor set and Euclidean space $\mathbb{R}^{\dim\Lambda}$. Moreover, Williams completely studied the dynamics of the restriction of a diffeomorphism on expanding attractor. In particular, he proved that this restriction is conjugate to the shift map of a generalized solenoid, see surveys \cite{GrinesPochinkaZh2014,GrinesZh2006}. Interesting examples of expanding attractors was constructed by Farrel and Jones \cite{FarrelJones1980,FarrelJones1981,Jones1983,Jones1986}. They consider interior dynamics with no embedding of expanding attractors in supporting manifolds. Robinson and Williams \cite{RobinsonWilliams76} constructed two diffeomorphisms $f$ and $g$ of non-homeomorphic 5-manifolds with expanding 2-dimensional attrac\-tors $\Lambda _f$ and $\Lambda _g$ respectively such that the restriction $f|_{\Lambda_f}: \Lambda_f\to\Lambda_f$ is conjugate to the restriction
$g|_{\Lambda_g}: \Lambda_g\to\Lambda_g$ but there is not even a homeomorphism from a neighborhood of $\Lambda_f$ to a neighborhood of $\Lambda_g$ taking $\Lambda_f$ to $\Lambda_g$. Another examples see in \cite{Bothe1992,IsaenkovaZh2009}.

Similar notation holds for axiom A flows (in short, A-flows) \cite{HirchPalisPughShub70}. By definition, a basic set is nontrivial if it not an isolated trajectory. In particular, a nontrivial basic set is not a fixed point (singularity) of A-flow. Since trajectories of flow are one-dimensional, any nontrivial basic set of A-flow has the topological dimension no less than one, and a supporting manifold admitting a nontrivial basic set has the dimension no less than three. Moreover, a nontrivial one-dimensional basic set is of saddle type while one-dimensional attractors and repellers on a 3-manifold are always trivial. One can prove that a two-dimensional basic set on $M^3$ is either an attractor or repeller, and two-dimensional attractors and repellers are exactly expanding attractors and contracting repellers respectively \cite{MedvedevZh2022}.

Questions of embedding and classification of orientable codimension one expanding attractors for A-diffeomorphisms was completely investigated by Grines \cite{Gr75-77}, Grines and Zhuzhoma \cite{GinesZhuzhoma1979,GinesZhuzhoma2005}, and Plykin \cite{Plykin74,Plykin84}. As to another dimensions in frame of our knowledge, there are results concerning one-dimensional expanding attractors in 3-manifolds \cite{JingNiWang2004,JimingMaBinYu2007,JimingMaBinYu2011}.

In the paper, we construct structurally stable diffeomorphisms with expanding attractors of various dimension on $n$-sphere $\mathbb{S}^n$ and $n$-torus $\mathbb{T}^n$, $n\geq 2$. We also prove the existence of A-diffeomorphisms with few-dimensional expanding attractors on arbitrary closed manifolds. The main results concerning diffeomorphisms are the following statements.
Below, $[x]$ means the integer part of the number $x$.
\begin{theorem}\label{thm:one-to-antie-dim-on-sphere}
Let $\mathbb{S}^n$, $n\geq 2$, be an $n$-sphere. Then given any $d\in\{1,\ldots,[\frac{n}{2}]\}$, there is a structurally stable diffeomorphism $f: \mathbb{S}^n\to\mathbb{S}^n$ with a $d$-dimensional non-orientable expanding attractor.
\end{theorem}

\begin{theorem}\label{thm:one-to-cod-one-dim-on-torus}
Let $\mathbb{T}^n$, $n\geq 2$, be an $n$-torus. Then given any $1\leq d\leq n-1$, there is a structurally stable diffeomorphism $f: \mathbb{T}^n\to\mathbb{T}^n$ with a $d$-dimensional orientable expanding attractor.
\end{theorem}

\begin{theorem}\label{thm:one-to-antie-dim-on-any-manifold}
Given any closed $n$-manifold $M^n$, $n\geq 2$, and any $d\in\{1,\ldots,[\frac{n}{2}]\}$, there is an axiom A diffeomorphism $f: M^n\to M^n$ with a $d$-dimensional non-orientable expanding attractor.
\end{theorem}

Let us consider expanding attractors of A-flows. It follows from \cite{MedvedevZh2020-translation} for $n=3$ and \cite{MedvedevZh2022}, Theorem 4, for $n\geq 4$ that given any closed orientable $n$-manifold $M^n$, $n\geq 3$, there is an A-flow $f^t$ on $M^n$ such that the non-wandering set $NW(f^t)$ contains a two-dimensional attractor. As to $d$-dimensional expanding attractors for $d\geq 3$, we prove the following statements.
\begin{theorem}\label{thm:d-dim-for-a-flows}
Given any $n$-sphere $\mathbb{S}^n$, $n\geq 2d+2\geq 8$, $d\geq 3$, there is an A-flow $\varphi^t$ on $\mathbb{S}^n$ such that the non-wandering set $NW(\varphi^t)$ contains a $d$-dimensional expanding attractor. In addition, $varphi^t$ has an isolated hyperbolic source.
\end{theorem}

Remark that it follows from Theorem \ref{thm:one-to-cod-one-dim-on-torus} that given any $2\leq d\leq n-1$, $n\geq 3$, there is a closed $n$-manifold $M^n$ supporting a structurally stable A-flow with an orientable $d$-dimensional expanding attractor. To be precise, $M^n$ is a mapping torus supporting a dynamical suspension over a structurally stable diffeomorphism $f: \mathbb{T}^{n-1}\to\mathbb{T}^{n-1}$ with a $(d-1)$-dimensional orientable expanding attractor.

\textsl{Acknowledgments}. This work is an output of a research project implemented as part of the Basic Research Program at the National Research University Higher School of Economics.

\section{Previous results and basic definitions}\label{s:prev}

\textsl{A-diffeomorphisms}.
A diffeomorphism $f: M^n\to M^n$ is called an \textit{A-diffeomorphism} if its non-wandering set $NW(f)$ is hyperbolic and periodic points are dense in $NW(f)$ \cite{Smale67}. The existence of hyperbolic structure implies the existence of unstable $W^u(x)$ and stable $W^s(x)$ manifolds respectively for every point $x\in NW(f)$ \cite{HirschPughShub77-book}.
Due to Smale's Spectral Decomposition Theorem, the non-wandering set $NW(f)$ is a finite union of pairwise disjoint $f$-invariant closed sets $\Omega_1$, $\ldots$, $\Omega_k$
such that every restriction $f|_{\Omega_i}$ is topologically transitive. These $\Omega_i$ are called the \textit{basic sets} of $f$. The dimension $\dim W^u(x)$, $x\in\Omega_i$ s called a \textit{Morse index} of $\Omega_i$.
A basic set $\Omega$ is an \textit{expanding attractor} if $\Omega$ is an attractor and the topological dimension $\dim\Omega$ equals Morse's index of $\Omega$ \cite{Williams1974}. A basic set $\Lambda$ of $f$ is called a \textit{contracting repeller} provided $\Lambda$ is an expanding attractor of $f^{-1}$.

By definition, let $W^s_{\epsilon}(x)\subset W^s(x)$ (resp. $W^u_{\epsilon}(x)\subset W^u(x))$ be the $\epsilon$-neighborhood of $x$ in the intrinsic topology of the manifold $W^s(x)$ (resp. $W^u(x))$, where $\epsilon >0$. We say that a basic set $\Omega$ is \textit{orientable} provided the index of intersection $W^s_{\alpha}(x)\cap W^u_{\beta}(x)$ is the same at each point of this intersection for any $\alpha>0$, $\beta>0$, $x\in\Omega$. Smale's solenoid is an orientable attractor \cite{Smale67} while Plykin attractor is a non-orientable one \cite{Plykin74}.

\medskip
\textsl{Structural stability and $\Omega$-stability}.
Let $\mathit{Diff}^1~(M^n)$ be the space of $C^1$ diffeomorphisms on $M^n$ endowed with the uniform $C^1$ topology \cite{Hirsch-book-1976}. Recall that diffeomorphisms $f$, $g\in Diff^1(M)$ are (topologically) \textit{conjugate} if there is a homeomorphism $\varphi : M\to M$ such that $\varphi \circ f = g\circ \varphi$. A diffeomorphism $f\in Diff^1(M)$ is called \textit{structurally stable} if there is a neighborhood $U(f)\subset Diff^1(M)$ of $f$ such that any $g\in U$ is conjugate to $f$.

Let $W_1$, $W_2\subset M^n$ be two immersed submanifolds. One says that $W_1$, $W_2$ are intersected \textit{transversally} provided given any point $x\in W_1\cap W_2$, the tangent bundles $T_xW_1$, $T_xW_2$ generate the tangent bundle $T_xM^n$. In this case $\dim T_xW_1 + \dim T_xW_2\ge \dim T_xM^n$. According Mane \cite{Mane88} and Robinson \cite{Robinson76}, an A-diffeomorphism $f$ is structurally stable if and only if invariant manifolds $W^s(x)$, $W^u(y)$ are intersected transversally for any $x$, $y\in NW(f)$. The last condition is called a \textit{strong transversality condition}.

\medskip
\textsl{Morse-Smale and DA-diffeomorphisms}.
In 1960, Steve Smale \cite{Smale60a} introduced a class of dynami\-cal systems (flows and diffeomorphisms) called later Morse-Smale systems. By definition, a diffeomorphism $f$ is Morse-Smale provided the non-wandering set $NW(f)$ consists of a finitely many periodic orbits, each periodic orbit is hyperbolic and, stable and unstable manifolds of periodic orbits intersect transversally.
One can define Morse-Smale systems as being those that are structurally stable and have non-wandering sets that consist of a finite number of orbits.
Smale \cite{Smale60a} discovered a deep connections between dynamics and the topological structure of support manifolds, see also \cite{Smale61a}. In particular, he proved the following statement which we need later on.
\begin{lemma}\label{lm:ms-grad-like-exist}
Any closed manifold $M^n$ admits a gradient-like Morse-Smale diffeomorphism with an isolated sink.
\end{lemma}

A DA-diffeomorphism $g: \mathbb{T}^n\to\mathbb{T}^n$ is an A-diffeomorphism such that either $NW(g)$ consists of codimension one orientable expanding attractor and finitely many isolated source periodic orbits or $NW(g)$ consists of codimension one orientable contracting repeller and finitely many isolated sink periodic orbits. Mainly, we'll consider a DA-diffeomorphism with a codimension one orientable expanding attractor called a DA-attractor. Such simplest DA-diffeomorphism can be obtained by Smale's surgery from a codimension one Anosov diffeomorphisms $A: \mathbb{T}^n\to\mathbb{T}^n$ such that $g_*=A_*: H_1(\mathbb{T}^n)\to H_1(\mathbb{T}^n)$, see details in \cite{Robinson-book-99}. In this case, $g$ has a DA-attractor $\Lambda_a$ that equal to $\mathbb{T}^n\setminus W^u(p)$ where $W^u(p)$ is the unstable manifold of a unique isolated source $p$. Later on, one needs the following statement proved in general case by Plykin \cite{Plykin84}.
\begin{lemma}\label{lm:cod-one-da-exist}
Given any $n\geq 2$, there is a structurally stable DA-diffeomorphism $f: \mathbb{T}^n\to\mathbb{T}^n$ with an $(n-1)$-dimensional orientable expanding attractor.
\end{lemma}

\medskip
Given maps $f_i: M_i\to M_i$, $i=1,2,\ldots,m$, we denote by
 $$F_{12\ldots m}=(f_1,f_2,\ldots,f_m): M_1\times M_2\times\cdots\times M_m\to M_1\times M_2\times\cdots\times M_m$$
the map
 $$F_{12\ldots m}(x_1,x_2,\ldots,x_m)=(f_1(x_1),f_2(x_2),\ldots,f_m(x_m)),\quad x_i\in M_i, \,\, i=1,\ldots,m.$$

\begin{lemma}\label{lm:previous}
Let $f_i: M^{k_i}_i\to M^{k_i}_i$ be a structurally stable diffeomorphism (resp., A-diffeomorp\-hism) of closed $k_i$-manifold $M^{k_i}$, $k_i\geq 1$, with non-wandering set $NW(f_i)$ consisting of basic sets $\Omega^{(i)}_1$, $\ldots$, $\Omega^{(i)}_{l_i}$, $i=1,2,\ldots,m$. Then the mapping
 $$F_{12\ldots m}: M_1^{k_1}\times M_2^{k_2}\times\cdots\times M_m^{k_m}\to M_1^{k_1}\times M_2^{k_2}\times\cdots\times M_m^{k_m}$$
is a structurally stable (resp., A-diffeomorphism) with non-wandering set
 $$NW(F_{12\ldots m})=NW(f_1)\times NW(f_2)\times\cdots\times NW(f_m).$$
Moreover, the spectral decomposition of $F_{12\ldots m}$ consists of basic sets $\Omega^{(1)}_{i_1}\times\Omega^{(2)}_{i_2}\times\cdots\times\Omega^{(m)}_{i_m}$ where $1\leq i_s\leq l_s$, $s=1,\ldots,m$
\end{lemma}
\textsl{Proof}. It is enough to prove the statement for $m=2$. Recall that a non-wandering set is the topological closure of set of periodic points. Therefore given any point $(x,y)\in NW(f_1)\times NW(f_2)$, there is $(p_1,p_2)\in NW(f_1)\times NW(f_2)$ arbitrary close to $(x,y)$ where $p_i$ is a periodic point of $f_i$, $i=1,2$. Obviously, $(p_1,p_2)$ is a periodic point of $F_{12}$. Hence, $NW(f_1)\times NW(f_2)\subset NW(F_{12})$. Due to \cite{Wiseman2017}, $NW(F_{12})\subset NW(f_1)\times NW(f_2)$. As a consequence, $NW(F_{12})=NW(f_1)\times NW(f_2)$.

The stable and unstable bundles of $f_1$ and $f_2$ form the stable and unstable bundles of $F_{12}$ as follows $\mathbb{E}_{F_{12}}^s=\mathbb{E}_{f_1}^s\oplus\mathbb{E}_{f_2}^s$, $\mathbb{E}_{F_{12}}^u=\mathbb{E}_{f_1}^u\oplus\mathbb{E}_{f_2}^u$. Since the differential $DF_{12}$ has the block structure $\left(D(f_1),D(f_2)\right)$, the fiber bundle decomposition $\mathbb{E}_{F_{12}}^s\oplus\mathbb{E}_{F_{12}}^u$ is invariant under $DF_{12}$, and $NW(F_{12})=NW(f_1)\times NW(f_2)$ has hyperbolic structure.

Let $\Omega_i$ be a basic set of $f_i$, $i=1,2$. Clearly, $\Omega_1\times\Omega_2=\Omega_{12}$ is a closed invariant set of $F_{12}$. We have to prove that $\Omega_{12}$ is a transitive set of the restriction $F_{12}|_{\Omega_{12}}$. First, suppose that the both $\Omega_1$ and $\Omega_2$ are trivial, i.e. $\Omega_1$ and $\Omega_2$ are periodic orbits. Then $\Omega_{12}$ is a periodic point, i.e. $\Omega_{12}$ is a basic set. Suppose now that the both $\Omega_1$ and $\Omega_2$ are nontrivial. Let $V_1$, $V_2$ be relatively open sets in $\Omega_{12}$. Without loss of generality, one can assume that $V_1$ is a rectangle, i.e. $V_1=V_{1x}\times V_{1y}$ where $V_{1x}$ and $V_{1y}$ are relatively open sets in $\Omega_1$ and $\Omega_2$ respectively. Similarly, $V_2=V_{2x}\times V_{2y}$. According \cite{Anos70,Bow71}, some iteration of A-diffeomorphism on a nontrivial basic set is topologically mixing. Since a transitivity under some iteration implies a transitivity, one can assume that the restriction $f_i|_{\Omega_i}$, $i=1,2$, is a topologically mixing mapping. Therefore there is $n_1$ such that $f^m_1(V_{1x})\cap V_{2x}\neq\emptyset$ for all $m\geq n_1$. Similarly, there is $n_2$ such that $f^m_2(V_{1y})\cap V_{2y}\neq\emptyset$ for all $m\geq n_2$. This follows that $F_{12}^{m}(V_1)\cap V_2\neq\emptyset$ for all $m\geq\max\{n_1,n_2\}$. We see that $F_{12}|_{\Omega_{12}}$ is topologically mixing, and as consequence $F_{12}|_{\Omega_{12}}$ is transitive (actually, one proves that $\Omega_1\times\Omega_2=\Omega_{12}$ is a nontrivial basic set). When $\Omega_1$ is a nontrivial basic set while $\Omega_2$ is a trivial basic set, one can prove similarly that $\Omega_1\times\Omega_2=\Omega_{12}$ is a nontrivial basic set as well.

Let us prove that $F_{12}$ is structurally stable provided $f_1$ and $f_2$ are structurally stable. It is enough to check a strong transversality condition.. Suppose an unstable manifold $W^u_{F_{12}}(p_1,q_1)$ intersects a stable manifold $W^s_{F_{12}}(p_2,q_2)$ at some point $(x,y)\in M_1\times M_2$. Since $f_i: M_i\to M_i$, $i=1,2$, is structurally stable, the unstable manifold $W^u_{f_1}(p_1)$ intersects transversally the stable manifold $W^s_{f_1}(p_2)$ at $x$ and the unstable manifold $W^u_{f_2}(q_1)$ intersects transversally the stable manifold $W^s_{f_2}(q_2)$ at $y$. It follows from $W^u_{F_{12}}(p_1,q_1)=W^u_{f_1}(p_1)\times W^u_{f_2}(q_1)$ and $W^s_{F_{12}}(p_2,q_2)=W^s_{f_1}(p_2)\times W^s_{f_2}(q_2)$ that $W^u_{F_{12}}(p_1,q_1)$ intersects transversally $W^s_{F_{12}}(p_2,q_2)$ at $(x,y)$. Hence, $F_{12}$ is structurally stable.
$\Box$

\begin{lemma}\label{lm:most-general}
Let $f_i: M^{k_i}_i\to M^{k_i}_i$ be a structurally stable (resp., A-diffeomorphism) of closed $k_i$-manifold $M^{k_i}$, $k_i\geq 1$, with a $d_i$-dimensional expanding attractor $\Lambda_i$, $i=1,2,\ldots,m$. Then the structurally stable diffeomorphism (resp., A-diffeomorphism)
 $$F_{12\ldots m}: M_1^{k_1}\times M_2^{k_2}\times\cdots\times M_m^{k_m}\to M_1^{k_1}\times M_2^{k_2}\times\cdots\times M_m^{k_m}$$
has the $(d_1+\cdots +d_m)$-dimensional expanding attractor $\Lambda_1\times\cdots\times\Lambda_m$. In addition, $F_{12\ldots m}$ has an isolated source.
\end{lemma}
\textsl{Proof}. Due to Lemma \ref{lm:previous}, $\Lambda_1\times\cdots\times\Lambda_m$ is a basic set of $F_{12\ldots m}$. Since $\Lambda_i$ is an attractor for any $1\leq i\leq m$, $\Lambda_1\times\cdots\times\Lambda_m$ is an attractor. Since $\Lambda_i$ is expanding, $\dim\Lambda_i=\dim W^u_{f_i}(x_i)$, $x_i\in\Lambda_i$. It follows from Lemma \ref{lm:previous} that $W^u_F(x_1,\ldots,x_m)=W^u_{f_1}(x_1)\times\cdots\times W^u_{f_m}(x_m)$.
As a consequence,
 $$ \dim(\Lambda_1\times\cdots\times\Lambda_m)=d_1+\cdots +d_m=\dim W^u_{f_1}(x_1)+\cdots +\dim W^u_{f_m}(x_m)=\dim W^u_F(x_1,\ldots,x_m). $$
Hence, $\Lambda_1\times\cdots\times\Lambda_m$ is an expanding attractor of dimension $(d_1+\cdots +d_m)$.
$\Box$

\begin{corollary}\label{cor:with-sink}
Suppose the condition of Lemma \ref{lm:most-general} holds. Let $f_{m+1}: M^{k_{m+1}}\to M^{k_{m+1}}$ be a structurally stable diffeomorphism with an isolated sink $\Omega$. Then $$\Lambda_1\times\cdots\times\Lambda_m\times\{\omega\}\subset M_1^{k_1}\times M_2^{k_2}\times\cdots\times M_m^{k_m}\times M_{m+1}^{k_{m+1}}$$
is $(d_1+\cdots +d_m)$-dimensional expanding attractor of $F_{12\ldots m,m+1}$.
\end{corollary}

\section{Proofs of main results}\label{s:proof}

\medskip
\textsl{Proof of Theorem \ref{thm:one-to-antie-dim-on-sphere}}. First, we prove the following statement concerning one-dimensional expanding attractors.
\begin{proposition}\label{prop:one-dim-on-any-sphere}
Any $n$-sphere $\mathbb{S}^n$, $n\geq 2$, admits a structurally stable diffeomorphism $\mathbb{S}^n\to\mathbb{S}^n$ with non-orientable one-dimensional expanding attractor.
\end{proposition}
\textsl{Proof of Proposition \ref{prop:one-dim-on-any-sphere}}. Well-known that there is a structurally stable diffeomorphism $f: \mathbb{S}^2\to\mathbb{S}^2$ with Plykin's attractor $\Lambda_a$ which is a non-orientable one-dimensional expanding attractor \cite{Plykin74}. Thus for $n=2$ the statement is true. Consider $\mathbb{S}^2$ to be the diameter of $\mathbb{S}^3$. Then $\mathbb{S}^2$ divides $\mathbb{S}^3$ into closed 3-disks $D_1$, $D_2$ such that $\mathbb{S}^2=D_1\cap D_2$. It is convenient to consider each $D_1$, $D_2$ to be the unit disk in $\mathbb{R}^3$ endowed with the polar coordinates $(\rho,\phi,\theta)$ where $\mathbb{S}^2$ is identified with $\rho=1$. The differential equation $\dot{\rho}=\rho(1-\rho)$ at each ray $\phi=\phi_0$, $\theta=\theta_0$ defines a flow $g^t$ on $D_1\cup D_2=\mathbb{S}^3$. Note that the set of fixed points of $g^t$ consists of two sources
$\alpha_1\in D_1$, $\alpha_2\in D_2$ and the sphere $\mathbb{S}^2$. Let $g=g^1$ be the shift along the trajectories on the time $t=1$. The calculations show that $\rho'(0)=e>1$ and $\rho'(1)=\frac{1}{e}<1$. Therefore, $\alpha_1$, $\alpha_2$ are hyperbolic sources. The sphere $\mathbb{S}^2$ is an attractive set of $g$, and any point of $\mathbb{S}^2$ is a hyperbolic sink in the $\rho$-direction.

There is a tubular neighborhood $T(\mathbb{S}^2)\subset\mathbb{S}^3$ of $\mathbb{S}^2$ which is trivial fiber bundle $\mathbb{S}^2\times [-1;1]$ with the base $\mathbb{S}^2\times\{0\}=\mathbb{S}^2\subset\mathbb{S}^3$ and fibers are arcs of trajectories of $g^t$. By construction, Plykin's diffeomorphisms $f: \mathbb{S}^2\to\mathbb{S}^2$ preserves orientation. Therefore there is a diffeotopy $f_{\alpha}: \mathbb{S}^2\to\mathbb{S}^2$, $0\leq\alpha\leq 1$, such that $f_0=f$, $f_1=id$. Let $\theta: \mathbb{R}\to [0;1]$ be a smooth even function such that $\theta(0)=0$ and $\theta(x)=1$ for $x\geq 1$. Define the diffeomorphism $F: \mathbb{S}^3\to\mathbb{S}^3$ as follows. For $\mathbb{S}^3\setminus T(\mathbb{S}^2)$, set $F=g$. For $(x;\beta)\in\mathbb{S}^2\times [-1;1]=T(\mathbb{S}^2)$, put by definition,
  $$F: T(\mathbb{S}^2)=\mathbb{S}^2\times [-1;1]\to\mathbb{S}^2\times [-1;1],\quad F(x,\beta)=(f_{\alpha(\theta(\beta))}(x),g(\beta)).$$
We see that $F(x,0)=(f(x),0)$. Since $f_{\alpha(\theta(\pm 1))}=id$, $F|_{\partial T(\mathbb{S}^2}=(id,g(\cdot))$. Hence, $F: \mathbb{S}^3\to\mathbb{S}^3$ is a structurally stable diffeomorphism with the non-orientable one-dimensional expanding attractor $\Lambda_a$. Similarly, one can construct a structurally stable diffeomorphism $\mathbb{S}^4\to\mathbb{S}^4$ with the non-orientable one-dimensional expanding attractor $\Lambda_a$. Continuing by this way, one can get a diffeomorphism desired for any $n\geq 3$. This completes the proof of Proposition \ref{prop:one-dim-on-any-sphere}.
$\diamondsuit$

After Proposition \ref{prop:one-dim-on-any-sphere}, we see that the result is true for $n=2,3$. Due to Lemma \ref{lm:most-general} and Proposition \ref{prop:one-dim-on-any-sphere}, the manifold $\mathbb{S}^2\times\mathbb{S}^i$ admits a structurally stable diffeomorphism $F_{1i}=(f_1,f_i): \mathbb{S}^2\times\mathbb{S}^i\to\mathbb{S}^2\times\mathbb{S}^i$ with a two-dimensional expanding attractor $\Lambda_i^{(2)}=\Lambda\times\Lambda_i^{(1)}$ for any $i\geq 2$ where $\Lambda$ is a one-dimensional (Plykin's) expanding attractor of $f_1: \mathbb{S}^2\to\mathbb{S}^2$ while $\Lambda_i^{(2)}$ is a one-dimensional attractor of $f_i: \mathbb{S}^i\to\mathbb{S}^i$. It follows from Lemma \ref{lm:most-general} that $\mathbb{S}^2$ (resp., $\Lambda_i^{(2)}$) has an isolated source $\alpha$ of $f_1$ (resp., $\alpha_i$ of $f_i$). Therefore, there is a closed 2-disk $D^2\subset\mathbb{S}^2$ containing $\Lambda$ such that $f_1(D^2)\subset D^2$, $\alpha\notin D^2$. Without loss of generality, one can assume that $D^2$ has a smooth boundary. Similarly, there is a closed $i$-disk $D^i\subset\mathbb{S}^i$ with a smooth boundary such that $f_i(D^i)\subset D^i$, $\Lambda_i^{(2)}\subset D^i$, and $\alpha_i\notin D^i$. Hence, $F_{1i}(D^2\times D^i)\subset D^2\times D^i$ and $\Lambda_i^{(2)}\subset D^2\times D^i$. The $(2+i)$-disk $D^2\times D^i$ has a smooth boundary. Gluing a $(2+i)$-disk $D^{2+i}$ along its boundary to $D^2\times D^i$, one gets $\mathbb{S}^{2+i}$. Since $F_{1i}(D^2\times D^i)\subset D^2\times D^i$, we can extend $F_{1i}$ to $\mathbb{S}^{2+i}$ to get a structurally stable diffeomorphism $F_i: \mathbb{S}^{2+i}\to\mathbb{S}^{2+i}$ with an isolated source in $D^{2+i}$. By construction, $F_i$ has the two-dimensional expanding attractor $\Lambda_i^{(2)}$, $i\geq 2$.
We see that every $\mathbb{S}^n$, $n\geq 4$, admits a structurally stable diffeomorphism with a two-dimensional expanding attractor. Thus for $n=4,5$, the result is true. This follows that $\mathbb{S}^2\times\mathbb{S}^i$ admits a structurally stable diffeomorphism with a three-dimensional expanding attractor for any $i\geq 4$. Similarly, one can prove that $\mathbb{S}^n$, $n\geq 5$, admits a structurally stable diffeomorphism with a three-dimensional expanding attractor (so, the result is true for $n=6,7$). Continuing by this way, one get the completes proof.
$\Box$

Remark. By construction, $f$ has an isolated source.

\medskip
\textsl{Proof of Theorem \ref{thm:one-to-cod-one-dim-on-torus}}. If $d=n-1$, then $f$ is a DA-diffeomorphism, see Lemma \ref{lm:cod-one-da-exist}. Suppose $d+1<n$ and $n\geq 3$. Let $f_1: \mathbb{T}^{d+1}\to\mathbb{T}^{d+1}$ be a DA-diffeomorphism with $d$-dimensional orientable expanding attractor $\Lambda_a$. Due to Lemma \ref{lm:ms-grad-like-exist}, there is a gradient-like Morse-Smale diffeomorphism $f_2: \mathbb{T}^{n-d-1}\to\mathbb{T}^{n-d-1}$ with a sink, say $\omega$. It follows from Corollary \ref{cor:with-sink} that
 $$F_{12}=(f_1,f_2): \mathbb{T}^{d+1}\times\mathbb{T}^{n-d-1}=\mathbb{T}^n\to\mathbb{T}^{d+1}\times\mathbb{T}^{n-d-1}=\mathbb{T}^n$$
is a structurally stable diffeomorphism with the $d$-dimensional orientable expanding attractor $\Lambda\times\{\omega\}$.
$\Box$

\medskip
\textsl{Proof of Theorem \ref{thm:one-to-antie-dim-on-any-manifold}}. Take a closed $n$-manifold $M^n$, $n\geq 2$. According Lemma \ref{lm:ms-grad-like-exist}, $M^n$ admits a gradient-like Morse-Smale diffeomorphism
$\psi: M^n\to M^n$ with an isolated sink $\omega$. Fix some $d\in\{1,\ldots,[\frac{n}{2}]\}$. Due to Theorem \ref{thm:one-to-antie-dim-on-sphere}, $n$-sphere $\mathbb{S}^n$ admits a structurally stable diffeomorphism
$f: \mathbb{S}^n\to\mathbb{S}^n$ with a $d$-dimensional non-orientable expanding attractor $\Lambda_a$. By construction, $f$ has an isolated source $\alpha$ (see Remark 1). Let $B_{\omega}\subset M^n$ be an $n$-ball containing $\Omega$ such that $\psi(B_{\omega})\subset B_{\omega}$. Let $B_{\alpha}\subset\mathbb{S}^n$ be an $n$-ball containing $\alpha$ such that $f^{-1}(B_{\alpha})\subset B_{\alpha}$ and $B_{\alpha}$ does not contain any non-wandering points of $f$ except $\alpha$. Let $M^n\sharp\mathbb{S}^n$ be a connected sum obtained from $M^n\setminus B_{\omega}$ and $\mathbb{S}^n\setminus B_{\alpha}$ after identifying along the boundary components $\partial B_{\omega}$, $\partial B_{\alpha}$. Since $\omega$ is a sink and $\alpha$ is a source, $\psi$ and $f$ induce a A-diffeomorphism $f: M^n\sharp\mathbb{S}^n\to M^n\sharp\mathbb{S}^n$ with the $d$-dimensional non-orientable expanding attractor $\Lambda_a$. Since $M^n\sharp\mathbb{S}^n=M^n$, $f$ is a diffeomorphism desired.
$\Box$

\medskip
\textsl{Proof of Theorem \ref{thm:d-dim-for-a-flows}}.
According Theorem \ref{thm:one-to-cod-one-dim-on-torus}, there is a structurally stable diffeomorphism $f: \mathbb{T}^d\to\mathbb{T}^d$ with a $(d-1)$-dimensional orientable expanding attractor. Hence, the dynamical suspension $sus^t(f)$ over $f$ is a structurally stable flow with a $d$-dimensional expanding attractor, say $\lambda_a$. The supporting manifold for $sus^t(f)$ is a mapping torus denoted by $M^{d+1}$ which is a closed $(d+1)$-manifold. Due to Whitney's Theorem \cite{Whitney1944}, we can assume that $M^{d+1}$ is smoothly embedded in $\mathbb{S}^{2d+2}$. A tubular neighborhood $T(M^{d+1})$ of $M^{d+1}$ is a locally trivial bundle with the base $M^{d+1}$ and fiber $\mathbb{B}^{d+1}$. One can assume that the boundary $\partial\left(T(M^{d+1})\right)$ of $T(M^{d+1})$ is a smooth submanifold. Therefore, there is a vector field $\vec v$ on $T(M^{d+1})$ which enters transversally through $\partial\left(T(M^{d+1})\right)$ in $T(M^{d+1})$ and coincides on $M^{d+1}$ with the vector field induced by the flow $sus^t(f)$. Moreover, one can assume that $M^{d+1}$ is an attracting invariant set of $\vec v$. Thus, $\lambda_a$ is an expanding attractor of $\vec v$. The Morse theory implies that there is an extension of $\vec v$ to a vector field $\vec V$ on $\mathbb{S}^{2d+2}$ such that $\vec V$ has only isolated hyperbolic fixed points including an isolated hyperbolic source. Due to Smale \cite{Smale61a}, one can assume that the restriction of $\vec V$ on $\mathbb{S}^{2d+2}\setminus T(M^{d+1})$ is a Morse-Smale vector field. Then $\vec V$ induces an A-flow $\varphi^t$ on $\mathbb{S}^{2d+2}$ with the $d$-dimensional expanding attractor $\lambda_a$. This proves Proposition for $n=2d+2$. Similarly to the proof of Theorem \ref{thm:one-to-antie-dim-on-sphere}, one can prove the statement desired for any $n\geq 2d+2$. This completes the proof.
$\Box$

\bigskip\noindent
\textit{E-mail:} medvedev-1942@mail.ru

\noindent
\textit{E-mail:} zhuzhoma@mail.ru

\end{document}